\documentclass[10pt,a4paper]{article}
\usepackage[utf8]{inputenc}
\usepackage{amsmath}
\usepackage{url}
\usepackage{amsfonts}
\usepackage{amssymb}
\usepackage{graphicx}
\usepackage[sectionbib,square]{natbib}
\usepackage[left=1in,right=1in,top=1in,bottom=1in]{geometry}
\author{Fernando Reynoso}
\title{Analytic Approach For Homogeneous\\and Non-Homogeneous Second Order PDE's\\with an Analytical Solution of\\ Navier Stokes Equations on $\mathbb{R}^3$\\ for a viscous incompressible fluid.}

\begin{document}
\maketitle
\begin{abstract}
This method take on consideration the Greens Identities and the Harmonic functions, allows the reduction of the order of the PDE(Partial differential equation) from second order case to a first order case, simplifying the analysis for linear and nonlinear PDE’s, not just giving an analytic solution for the behavior describe for the PDE, but at the same time give additional tools for existence and uniqueness for a broad range of PDE’s describe in theorem 1, testing this results in the Navier-Stokes Equation for free divergence fluid as well on other diffusive system as example, resulting on analytically well pose solution for each example.
\end{abstract}

\section{Description of the \textbf{Harmonic Reduction Method(HRM)}}

\hspace{0.20in}To fully understand this method is necessary to understand how harmonic functions behave and how they interact with the Laplace operator along a certain coordinate system, therefore, under this premise is necessary to define what is an harmonic function.

Definition: An harmonic function $\Gamma$ is a functions such that satisfy the Laplace differential equation, this means that fulfill the following condition.

\begin{equation*}
\Delta\Gamma=0
\end{equation*}

Notice that under this statement the harmonic function has multiple solutions that depend on the border conditions of the problems to define it's uniqueness, also notice that this type of function belong to a linear space, but the main goal of this paper is not to rediscover  how harmonics behave but use this results as a tool for other PDE's and find analytic solutions of them.

For this reason because harmonic functions belong to a linear space over the a certain domain with spatial proprieties, associated to constant coefficients related to the border conditions, to develop our analysis this constants will be replace by temporal functions, this way embedding this harmonic functions with a behavior over time, to describe behavior over time dependent PDE's.

Now as a second place, but no least important is good to review the Greens Identities, that are explicit consequence of the divergence theorem, therefore recalling this results are the following.

\begin{equation}
\int_Su\nabla v\cdot\hat{n}dS=\int_V\left(\nabla u\cdot\nabla v+ u\Delta v\right)dV\label{1Green}
\end{equation}

By taking in consideration the symmetry on how the gradient operator behave does imply in the Second Green Identity.

\begin{equation}
\int_Su\nabla v-v\nabla u\cdot\hat{n}dS=\int_V\left(u\Delta v-v\Delta u\right)dV\label{2Green}
\end{equation}

And.

\begin{equation}
\int_S\nabla\left(uv\right)\cdot\hat{n}dS=\int_V\Delta \left(uv\right) dV\label{Div-uv}
\end{equation}

base elements on the construction of this method.

Before continuing with the construction of this method,  is necessary to define the broad range of functions this method could be applied.

Let $\mathcal{F}$ be an homogeneous PDE of order two just as is described in the following equation.

\begin{align}
\mathcal{F}(t,\vec{x},u)=\Delta u+\vec{A}(t,\vec{x},u)\cdot\nabla u+B(t,\vec{x},u)\frac{\partial u}{\partial t}+C(t,\vec{x},u)u=0\label{GPDE}
\end{align}

By taking this considerations is implicit that the following results holds for any $v=\Gamma$ harmonic.

\begin{equation}
\int_S(u\nabla \Gamma)\cdot\hat{n}dS=\int_V\left(\nabla u\cdot\nabla\Gamma\right) dV
\end{equation}
\begin{equation}
\int_S(u\nabla \Gamma-\Gamma\nabla u)\cdot\hat{n}dS=-\int_V\Gamma\Delta udV
\end{equation}
\begin{equation}
\int_S(u\nabla \Gamma+\Gamma\nabla u)\cdot\hat{n}dS=\int_V\Delta(u\Gamma)dV
\end{equation}

Therefore by adding equation (6) and (7) is obtained the following result.

\begin{equation*}
\int_V2(\nabla u\cdot\nabla\Gamma) dV=\int_V\left[\Delta(u\Gamma)-\Gamma\Delta u\right]dV
\end{equation*}

Notice thanks to equation (\ref{GPDE}) by evaluating $\mathcal{F}(t,\vec{x},u\Gamma)$ our last result can be written as follows.

\begin{equation}
\int_V2\nabla u\cdot\nabla\Gamma dV=\int_V\left(\Gamma\vec{A}\cdot\nabla u+\Gamma B\partial_tu+\Gamma Cu-\vec{A_\Gamma}\cdot\nabla(u\Gamma)-B_\Gamma\partial_t(u\Gamma)-C_\Gamma(u\Gamma)\right)dV
\end{equation}

Notice that the notation was contracted by the use of the sub index $\Gamma$ referring to a change of that specific coefficient under the influence of $u$, this is done to take in consideration the definition of $\mathcal{F}$ described at equation (\ref{GPDE}), therefore by applying some algebra is found the following result.

\begin{equation*}
\int_V\left[\Gamma\left((\vec{A}-\vec{A_\Gamma})\cdot\nabla u+(B-B_\Gamma)\partial_tu+(C-C_\Gamma)u\right)-u\left(\vec{A_\Gamma}\cdot\nabla\Gamma+B_\Gamma\partial_t\Gamma\right)-2\nabla u\cdot\nabla\Gamma\right]dV=0
\end{equation*}

As a consequence that is chosen any volume on the metric-space, make implicit that the expression to be integrated has to be 0, therefore implying the following.

\begin{equation}
2\nabla u\cdot\nabla\Gamma=\Gamma\left((\vec{A}-\vec{A_\Gamma})\cdot\nabla u+(B-B_\Gamma)\partial_tu+(C-C_\Gamma)u\right)-u\left(\vec{A_\Gamma}\cdot\nabla\Gamma+B_\Gamma\partial_t\Gamma\right)\label{reduce}
\end{equation}

From equation (\ref{reduce}) is important to notice that the function $\Gamma$ was chosen, therefore is a known parameter, as well shows that the PDE it self has reduce it's order from two to one by this procedure, facilitating not just the possibility to find a solution for the particular PDE, but also alows an extension to find proof for analytic and uniqueness analysis for certain PDE's, regardless if are linear or non linear PDE's.

\subsection*{Theorem 1: Harmonic Reduction for diffusive PDE's}
Any PDE that fulfills the following form. 

\begin{equation*}
\Delta u+\vec{A}(t,\vec{x},u)\cdot\nabla u+B(t,\vec{x},u)\partial_tu+C(t,\vec{x},u(t,\vec{x}))u=0
\end{equation*}

Can be reduce from a second order PDE as a system of PDE's related to a set of harmonic functions define in the domain of the PDE.

Proof: 

By replying the process shown previously, would be necessary for an solution candidate function $u$ define in a $n$-dimensional domain $\Omega$, therefore by using the result obtained at equation (\ref{reduce}).

\begin{equation*}
\Gamma_i\left((\vec{A}-\vec{A_{\Gamma_i}})\cdot\nabla u+(B-B_{\Gamma_i})\partial_tu+(C-C_{\Gamma_i})u\right)-u\left(\vec{A_{\Gamma_i}}\cdot\nabla\Gamma_i+B_{\Gamma_i}\partial_t\Gamma_i\right)=2\nabla u\cdot\nabla\Gamma_i
\end{equation*}

By rearrange the reduce PDE, after some algebra.

\begin{equation}
-u\left[\vec{A_{\Gamma_i}}\cdot\nabla\Gamma_i+B_{\Gamma_i}\partial_t\Gamma_i+\Gamma(C_{\Gamma_i}-C)\right]=\Gamma(B_{\Gamma_i}-B)\partial_tu+\nabla u\cdot\left[2\nabla\Gamma_i+\Gamma(A_{\Gamma_i}-A)\right]\label{Constrained-Reduce}
\end{equation}

By choosing $\Gamma_i=x_if_i(t)+1$ where $x_i$ is a coordinate along the $\hat{e_i}$ vector basis of the domain $\Omega$ and $f_i:[0,\infty)\rightarrow\mathbb{R}$ is a temporal function specific for that coordinate, implying after some algebra.

\begin{equation*}
u\left[\frac{(\vec{A_{\Gamma_i}}\cdot\hat{e_i})f_i(t)}{x_if_i(t)+1}+\frac{B_{\Gamma_i}x_i\partial_tf_i(t)}{x_if_i(t)+1}+(C_{\Gamma_i}-C)\right]
=(B-B_{\Gamma_i})\frac{\partial u}{\partial t}+\nabla u\cdot\left[(A-A_{\Gamma_i})-\frac{2f_i(t)\hat{e_i}}{x_if_i(t)+1}\right]
\end{equation*}

Notice that this expression is valid for any integer $i$ referencing the coordinate system of the studied domain, proving this way the theorem.

Is important to notice that in some cases the metric-space can be distorted in relation to a temporal function, as occurs in special and general relativity, for that reason the temporal partial derivative wasn't develop.

This differ from the usual reduction method because you don't know any prior solution, but you assigned a particular function that can be tuned by the temporal function embedded onto it, dealing with a set of ODE's that are quiet easy to solve in most cases instead of a single PDE that is hard to analyze because for some particular reason.

Now will be shown an example referring how to use this method.

To do this will be used as reference the Heat Equation on $\mathbb{R}^n\times[0,\infty)$ .

\begin{equation}
\Delta u-\alpha\partial_tu=0\label{Heat}
\end{equation}

On this particular case is easy to notice that $\vec{A}=\vec{0}$, $C=0$ and $B=-\alpha$ doesn't have dependence on $u$, therefore by equation (\ref{reduce}) is achieve the following result.

\begin{equation}
2\nabla u\cdot\nabla\Gamma_i=\alpha u\partial_t\Gamma_i\label{Heat-Reduce}
\end{equation}

Under this successful reduction by choosing $\Gamma_i=1+(x_i-x_i^\circ)f_i(t)$ with $x_i$ a coordinate of the system of coordinate fixed along time, and $f_i:[0,\infty)\rightarrow\mathbb{R}$ a temporal function left to posterior analysis, with this information does imply the following.

\begin{equation*}
2f_i\partial_iu=\alpha u(x_i-x_i^\circ)\partial_tf_i
\end{equation*}

Equivalent to.

\begin{equation}
\partial_iu=\frac{\alpha u(x_i-x_i^\circ)\partial_tf_i}{2f_i}\label{Heat-dir}
\end{equation}

This valid expression has many consequences such as an expression for the gradient.

\begin{equation}
\nabla u=\frac{\alpha u}{2}\sum_{i=1}^n\frac{(x_i-x_i^\circ)\partial_tf_i}{f_i}\hat{e_i}\label{Heat-Grad}
\end{equation}

Therefore the Laplacian becomes.

\begin{equation*}
\Delta u=u\left[\frac{\alpha^2}{4}\sum_{i=1}^n\left(\frac{(x_i-x_i^\circ)\partial_tf_i}{f_i}\right)^2+\frac{\alpha}{2}\sum_{i=1}^n\frac{\partial_tf_i}{f_i}\right]
\end{equation*}

Therefore by the studied PDE describe in equation (\ref{Heat}) the temporal differential would be describe by the following.

\begin{equation}
\frac{\partial u}{\partial t}=u\left[\frac{\alpha}{4}\sum_{i=1}^n\left(\frac{(x_i-x_i^\circ)\partial_tf_i}{f_i}\right)^2+\frac{1}{2}\sum_{i=1}^n\frac{\partial_tf_i}{f_i}\right]\label{Heat-temp}
\end{equation}

From equation (\ref{Heat-Grad}) is also deduce by doing a line integral.

\begin{equation}
u(t,\vec{x})=u(t,\vec{x_0})\exp\left[\frac{\alpha}{4}\sum_{i=1}^n\frac{|x_i-x_i^\circ|^2\partial _tf_i}{f_i}\right]
\end{equation}

By replacing thislast result on equation (\ref{Heat-temp}) does hold.

\begin{equation}
\frac{\partial_tu(t,\vec{x}_0)}{u(t,\vec{x}_0)}+\frac{\alpha}{4}\sum_{i=1}^n|x_i-x_i^\circ|^2\left[\frac{\partial}{\partial t}\left(\frac{\partial_tf_i}{f_i}\right)-\left(\frac{\partial_tf_i}{f_i}\right)^2\right]=\frac{1}{2}\sum_{i=1}^n\frac{\partial_tf_i}{f_i}\label{constrained-Heat-temp}
\end{equation}

Notice from this last equation that $u(t,\vec{x}_0)$ has no dependence of the coordinate system but the source of the heat, therefore the coefficient of each element related to the coordinate system has to become 0, thus this condition implies that for all $i$ does hold.

\begin{equation}
\frac{\partial}{\partial t}\left(\frac{1}{f_i}\frac{\partial f_i}{\partial t}\right)-\left(\frac{1}{f_i}\frac{\partial f_i}{\partial t}\right)^2=0\label{ODE-Heat-System}
\end{equation}

therefore by applying this method of resolution results on transforming a second order PDE to a system of non-linear ODE's that depends on a number of equations related to the spatial dimension of the space $\Omega$, for this particular case being $n=\dim\mathbb{R}^n$. Also important to notice that this system of ODE's will vary depending on the original PDE.

After solving the ODE describe in equation (\ref{ODE-Heat-System}), the solution describe by this method for the Heat equation becomes by setting the initial conditions for $f_i(0)=1$ and $\frac{\partial f_i}{\partial t}(0)=-1$.

\begin{equation*}
u(t,\vec{x})=u^\circ(\vec{x}_0)\sqrt{\left(\frac{1}{t+1}\right)^n}\exp\left[-\frac{\alpha |\vec{x}-\vec{x}_0|^2}{4(t+1)}\right]
\end{equation*}

Where $u^\circ(\vec{x}_0)$ is the initial data to be set on the PDE problem. 

Doing this is posible because we define the harmonic to fit our requirements for that reason $f_i$ is analyzed almost at the end of the process so their behavior could be taking as natural as posible from the behavior of the studied PDE.

Notice from this method have compatibility with other methods of first order PDE's such as the characteristic method, taking in consideration the information that you have under your disposal or can generate equivalence between border conditions, that will be shown in the next section.

\section{Equivalence between types of border conditions}

\hspace{0.20in}Notice that for BVP you can build up an equivalence between the Dirichlet and the Neumann, Robin's, Mixed and Cauchy Problem, that comes from equation (\ref{reduce}) that remains valid for all elements of the domain, in particular of the border.

\subsection*{Equivalence Theorem for BVP in the analysis of PDE}
\hspace{0.20in}For a PDE of the form of equation (\ref{GPDE}) associated to a BVP such that $B$ does not depend on the function $u$, then does exist an equivalence between border conditions from the Dirichlet data to Von Neumann data and the other way around. Implying this way a respective construction of Robin boundary conditions, or Cauchy boundary conditions.

Proof: Starting from equation (\ref{reduce}) is clear the following.

\begin{align*}
u\left(-\vec{A_\Gamma}\cdot\nabla\Gamma-B\partial_t\Gamma\right)+\Gamma\left((\vec{A}-\vec{A_\Gamma})\cdot\nabla u+(C-C_\Gamma)u\right)=2\nabla u\cdot\nabla\Gamma
\end{align*}

Therefore by setting the harmonic such that it gradient points along the interest surface does become this problem on an algebraic problem by using the available data to construct the other type of data that will be valid for every element of the border of the domain.

From this result then you can build the Cauchy or the Robin Border conditions as better fit the analysis of the studied PDE.

To show some application under a nonlinear case behavior, will be studied the following PDE.

\begin{equation}
\alpha\partial_tu=\Delta u+\beta u(1-u)\label{Fisher}
\end{equation}

That correspond to the reaction-diffusion model known as the Fisher Equation, also known as Kolmogorov–Petrovsky–Piskunov equation, that has a wide rage of applicability on natural sciences, being one of them the use of population control, this will be analyzed in $\mathbb{R}^3$ for a single specie of bacteria, for this case to keep things simple will be use spherical coordinates, with the Dirichlet boundary condition $u(t,\vec{r})=0$ for all $\vec{r}\in\partial\mathbb{R}^3$ for any $t\geq 0$ and will be symmetrical over the angular component of of the domain.

Now as a consequence of the previous analysis done the spatial derivative becomes the following by using $\Gamma=1+f/r$.

\begin{equation}
\partial_ru=\frac{\alpha ur\partial_tf}{2f}-\frac{\beta u^2(f+r)}{2}
\label{Spatial-Fisher-Der}
\end{equation}

Thus this results does hold the following border conditions $\partial_ru=0$ by letting $\hat{n}=\hat{r}$ on the surface integral of the gradient of $u$ does imply thanks to the divergence theorem that.

\begin{align}
\Delta u=0\label{Laplace-Fisher}
\end{align}

And.

\begin{align}
\partial_tu=-\frac{\beta}{\alpha}(u^2-u)\label{Constrained-Fisher}
\end{align}

Implying that the gradient becomes as a consecuence of the equation (\ref{Laplace-Fisher}).

\begin{align*}
\nabla u=\nabla\times A(t,\vec{r})
\end{align*}

Therefore those results imply the following partial solution of the PDE.

\begin{align*}
u(t,\vec{r})=u(t,\vec{r_0})+\int_{\vec{r_0}}^{\vec{r}}\nabla\times A(t,\vec{s})\cdot d\vec{s}
\end{align*}

By using this result on equation (\ref{Constrained-Fisher}) does imply the following expression.

\begin{align*}
\partial_tu(t,\vec{r_0})+\partial_tP(t,\vec{r})=-\frac{\beta}{\alpha}\left(u(t,\vec{r_0})^2-u(t,\vec{r_0})\right)-\frac{\beta}{\alpha}\left(P(t,\vec{r})^2+P(t,\vec{r})(2u(t,\vec{r_0})-1)\right)
\end{align*}

Where $P(t,\vec{r})=\int_{\vec{r}_0}^{\vec{r}}\nabla\times\vec{A}(t,\vec{s})d\vec{s}$. Also note that this last expression is equivalent to.

\begin{align*}
\partial_tu(t,\vec{r_0})+\frac{\beta}{\alpha}\left(u(t,\vec{r_0})^2-u(t,\vec{r_0})\right)=-\partial_tP(t,\vec{r})-\frac{\beta}{\alpha}\left(P(t,\vec{r})^2+P(t,\vec{r})(2u(t,\vec{r_0})-1)\right)
\end{align*}

Notice that this last PDE has to hold for any $\vec{r}$ in the domain, therefore by fixing the location over the source at $\vec{r_0}$, does imply the following first order ODE in for the temporal behavior on the source of the PDE.

\begin{align}
\partial_tu(t,\vec{r_0})=-\frac{\beta}{\alpha}(u^2(t,\vec{r_0})-u(t,\vec{r_0}))\label{Static-temporal-PDE-Fisher}
\end{align}

Therefore implying for an arbitrary location on the space the following ODE.

\begin{align}
\partial_tP(t,\vec{r})=-\frac{\beta}{\alpha}\left(P(t,\vec{r})^2+P(t,\vec{r})(2u(t,\vec{r_0})-1)\right)\label{NonStatic-temporal-PDE-Fisher}
\end{align}

From Equation (\ref{Static-temporal-PDE-Fisher}) the solution for that ODE, as a logistic temporal ODE is conclude the following.

\begin{align}
u(t,\vec{r_0})=\frac{u_0(\vec{r_0})}{\exp\left\lbrace-\frac{\beta t}{\alpha}\right\rbrace+u_0(\vec{r_0})\left(\exp1-\left\lbrace-\frac{\beta t}{\alpha}\right\rbrace\right)}\label{Fisher-Source Solution}
\end{align}

Notice that dividing the equation (\ref{NonStatic-temporal-PDE-Fisher}) by $P^2$ does imply in the next expression.

\begin{align*}
\frac{\partial_tP(t,\vec{r})}{P(t,\vec{r})^2}=-\frac{\beta}{\alpha}\left(1+\frac{(2u(t,\vec{r_0})-1)}{P(t,\vec{r})}\right)
\end{align*}

equivalent to the following.

\begin{align*}
\partial_tP(t,\vec{r})^{-1}-\frac{\beta(2u(t,\vec{r_0})-1)}{\alpha}P(t,\vec{r})^{-1}=\frac{\beta}{\alpha}
\end{align*}

By multiplying by the respective integration factor does hold.

\begin{align*}
\partial_t\left(\exp\left\lbrace-\frac{2\beta}{\alpha}\int_0^tu(t,\vec{r_0})dt\right\rbrace\exp\left\lbrace\frac{\beta t}{\alpha}\right\rbrace P(t,\vec{r})^{-1}\right)=\frac{\beta}{\alpha}\exp\left\lbrace-\frac{2\beta}{\alpha}\int_0^tu(t,\vec{r_0})dt\right\rbrace\exp\left\lbrace\frac{\beta t}{\alpha}\right\rbrace
\end{align*}

By replacing (\ref{Fisher-Source Solution}) to this result would lead to the following.

\begin{align*}
\partial_t\left(\frac{\exp\left\lbrace\frac{\beta t}{\alpha}\right\rbrace}{\left(1+u_0(\vec{r_0})\left(\exp\left\lbrace\frac{\beta t}{\alpha}\right\rbrace-1\right)\right)^2P(t,\vec{r})}\right)=\frac{\beta\exp\left\lbrace\frac{\beta t}{\alpha}\right\rbrace}{\alpha\left(1+u_0(\vec{r_0})\left(\exp\left\lbrace\frac{\beta t}{\alpha}\right\rbrace-1\right)\right)^2}
\end{align*}

By integrating over the time variable and after some algebra does hold.

\begin{align*}
\frac{\exp\left\lbrace\frac{\beta t}{\alpha}\right\rbrace}{\left(1+u_0(\vec{r_0})\left(\exp\left\lbrace\frac{\beta t}{\alpha}\right\rbrace-1\right)\right)^2P(t,\vec{r})}=\frac{1}{P(0,\vec{r})\left(1+u_0(\vec{r_0})\left(\exp\left\lbrace\frac{\beta t}{\alpha}\right\rbrace-1\right)\right)}
\end{align*}

Letting $P_0(\vec{r})=P(0,\vec{r})$ this last expression is equivalent to the following.

\begin{align*}
P(t,\vec{r})=\frac{P_0(\vec{r})\exp\left\lbrace\frac{\beta t}{\alpha}\right\rbrace}{1+u_0(\vec{r_0})\left(\exp\left\lbrace\frac{\beta t}{\alpha}\right\rbrace-1\right)}
\end{align*}

Where $P_0$ is the initial conditions for P. This allows to conclude that the partial solution is.

\begin{align}
u(t,\vec{r})=\frac{(u_0(\vec{r_0})+P_0(\vec{r}))\exp\left\lbrace\frac{\beta t}{\alpha}\right\rbrace}{1+u_0(\vec{r_0})\left(\exp\left\lbrace\frac{\beta t}{\alpha}\right\rbrace-1\right)}\label{Partial-Solution Fisher}
\end{align}

Therefore notice that to find pur desire function remains necessary to describe $P_0(\vec{r})$, define as the initial state of the line integral of the potential vector, being a explicit way to describe the behavior of the specie allowing to forecast how endanger a specie could be in a specific time over a certain habitat or certain conditions.

Notice that this construction allows a more manageable procedure, and a more comprehensive results to use on natural sciences, such as biology, chemistry and physics.

\section{A Look to Existence and Uniqueness}

\hspace{0.20in}Now is important to show how this algorithm interact with the existence and uniqueness of an arbitrary Homogeneous PDE as were define at (\ref{GPDE}).

For this reason is Known the theorem of existence and uniqueness of Picard-Lindel\"of, thinking in the following conditions:

\begin{enumerate}
\item[(1)] For the differential would be use the total differential as a reference for the function yet to define in each case denoted by $\mathcal{G}(t,\vec{x},u(t,\vec{x}))$, otherwise could be use an specific directional derivative compatible with the constrain of the problems.
\item[(2)] The respective set in consideration would be a subset of the domain $\Omega$, such that contain the initial conditions describe in the initial data according to the available information display in the PDE problem.
\item[(3)] Finally as a last consideration would be use the HRM to find the gradient and the Laplacian to induce the proper conditions to use the Picard-Lindel\"of Theorem.
\item[(4)] This function can be constrained along that path as if it was a single variable function of the type $\vec{\gamma}(t)=(t,\vec{x}(t))$.
\end{enumerate}

Therefore before applying this theorem is necessary to proof the extension of this one under the conditions describe.

\subsection*{Existence and Uniqueness Theorem for diffusive PDE's\\Extension of the Picard Lindel\"of Theorem}

\hspace{0.20in}Let $u\in \mathcal{C}_1([0,\infty))\times\Omega)$ belong to the family of functions that solve equation (\ref{GPDE}), and the initial data of the problem $u(t_0,\vec{x}_0)=u^\circ(\vec{x}_0)\in\mathcal{C}_1(\Omega)$ and $\mathcal{G}\in\mathcal{C}_1(I\times u(I))$, where $(t_0,\vec{x}_0)\in I$ a compact subset of $[0,\infty))\times\Omega$ and the total derivative is define as.

\begin{equation}
D_Tu=\partial_tu+\nabla u\cdot D_t\vec{x}=\mathcal{G}(t,\vec{x},u(t,\vec{x}))\label{Total}
\end{equation} 

Then does exist a unique solution for the PDE.

Proof: Notice that as a consequence of $u\in\mathcal{C}_1([0,\infty)\times\Omega)$ does exist the total derivative, being compatible with the PDE describe in equation (\ref{GPDE}), on the other side let's analyze the equation (\ref{Total}) in consideration of the studied type of PDE along this paper.

\begin{equation*}
\partial_tu+\nabla u\cdot D_T\vec{x}=\mathcal{G}(t,\vec{x},u(t,\vec{x}))
\end{equation*}

Notice that from the HRM can be deduce each partial differential with that information the gradient and the Laplace are constructed in terms of $u$, This constructions can be develop in deepness, but by recalling the result expressed on equation (\ref{Constrained-Reduce}).

\begin{align*}
-u\left(\vec{A_{\Gamma}}\cdot\nabla\Gamma+B_{\Gamma}\partial_t\Gamma+\Gamma(C_{\Gamma}-C)\right)=\Gamma(B_{\Gamma}-B)\partial_tu+\nabla u\cdot\left(2\nabla\Gamma+\Gamma(A_{\Gamma}-A)\right)
\end{align*}

We have to take in consideration 2 possible scenarios, first scenario $B$ depends on the function $u$, second scenario $B$ does not depend on $u$

On the first scenario would be enough to divide by the coefficient of the temporal partial derivative of $u$ resulting on the following.

\begin{align*}
\partial_tu=-\nabla u\cdot\left[\frac{2\nabla\Gamma+\Gamma(A_{\Gamma}-A)}{\Gamma(B_{\Gamma}-B)}\right]-u\left[\frac{\vec{A_{\Gamma}}\cdot\nabla\Gamma+B_{\Gamma}\partial_t\Gamma+\Gamma(C_{\Gamma}-C)}{\Gamma(B_{\Gamma}-B)}\right]
\end{align*} 

Notice that from the definition of $\mathcal{G}(t,\vec{x},u(t,\vec{x}))$ does hold.

\begin{equation*}
\mathcal{G}(t,\vec{x},u(t,\vec{x}))=\nabla u\cdot\left(D_T\vec{x}-\frac{2\nabla\Gamma+\Gamma(A_{\Gamma}-A)}{\Gamma(B_{\Gamma}-B)}\right)-u\left(\frac{\vec{A_{\Gamma}}\cdot\nabla\Gamma+B_{\Gamma}\partial_t\Gamma+\Gamma(C_{\Gamma}-C)}{\Gamma(B_{\Gamma}-B)}\right)
\end{equation*}

Therefore by adjusting $\Gamma$, such that the next condition holds.

\begin{align*}
D_T\vec{x}=\frac{2\nabla\Gamma}{\Gamma(B_{\Gamma}-B)}+\frac{(A_{\Gamma}-A)}{(B_{\Gamma}-B)}
\end{align*}

Then forces that $\mathcal{G}(t,\vec{x},u(t,\vec{x}))$ becomes.

\begin{align}
\mathcal{G}(t,\vec{x},u(t,\vec{x}))=-u\left(\frac{\vec{A_{\Gamma}}}{2}\cdot\left(D_T\vec{x}-\frac{(A_{\Gamma}-A)}{(B_{\Gamma}-B)}\right)+\frac{B_{\Gamma}\partial_t\Gamma}{(B_{\Gamma}-B)\Gamma}+\frac{(C_{\Gamma}-C)}{(B_{\Gamma}-B)}\right)\label{Total-Scenario1}	
\end{align}

Notice that $D_t\vec{x}\in C^K([0,\infty))\times\Omega)$, where $K$ is the lowest order of continuity between $\vec{A}$ and $B$.

On the mean time for the second scenario because $B_\Gamma=B$ can't be used the same type of procedure, but we know that holds the following.

\begin{align*}
\nabla u\cdot\left[\nabla\Gamma+\frac{\Gamma(\vec{A}_{\Gamma}-\vec{A})}{2}\right]=-u\left[\frac{\vec{A}_{\Gamma}\cdot\nabla\Gamma+B\partial_t\Gamma+\Gamma(C_{\Gamma_i}-C)}{2}\right]
\end{align*}

By translating this result onto the Einstein notation is known the following and assuming that $\Gamma_j$ correspond to the Harmonic related to the $j-$th direction.

\begin{align*}
\partial_ju\left[\partial_j\Gamma_j+\frac{\Gamma_j((\vec{A}_{\Gamma_j})_j-\vec{A_j})}{2}\right]=-u\left[\frac{(\vec{A}_{\Gamma_j})_j\partial_j\Gamma_j+B\partial_t\Gamma_j+\Gamma_j(C_{\Gamma_j}-C)}{2}\right]
\end{align*}

Therefore equivalent to.

\begin{align*}
\partial_ju=-u\left[\frac{(\vec{A_{\Gamma_j})_j}\partial_j\Gamma_j+B\partial_t\Gamma_j+\Gamma_j(C_{\Gamma_j}-C)}{2\partial_j\Gamma_j+\Gamma_j((\vec{A_{\Gamma_j}})-\vec{A})\cdot\hat{e_j}}\right]
\end{align*}

Therefore the gradient has to be.

\begin{align*}
\nabla u=-u\sum\frac{(\vec{A}_{\Gamma_j})_j\partial_j\Gamma_j+B\partial_t\Gamma_j+\Gamma_j(C_{\Gamma_j}-C)}{2\partial_j\Gamma_j+\Gamma_j(\vec{A}_{\Gamma_j}-\vec{A})\cdot\hat{e_j}}\hat{e_j}
\end{align*}

Implying the behavior of the laplacian as follows.

\begin{align*}
\Delta u=u\sum\left[\left(\frac{(\vec{A_{\Gamma_j})_j}\partial_j\Gamma_j+B\partial_t\Gamma_j+\Gamma_j(C_{\Gamma_j}-C)}{2\partial_j\Gamma_j+\Gamma_j((\vec{A_{\Gamma}})_j-\vec{A})\cdot\hat{e_j}}\right)^2-\partial_j\left(\frac{(\vec{A_{\Gamma_j})_j}\partial_j\Gamma_j+B\partial_t\Gamma_j+\Gamma_j(C_{\Gamma_j}-C)}{2\partial_j\Gamma_j+\Gamma_j((\vec{A_{\Gamma}})_j-\vec{A})\cdot\hat{e_j}}\right)\right]
\end{align*}

Therefore are define all the elements to describe the temporal partial derivative by doing the substitution of the laplacian and the gradient on the original PDE.

Leaving everything in terms of the coefficients and the function $u$ therefore ensuring the existence of $\mathcal{G}\in\mathcal{C}_0(\Omega\times[0,\infty))$, so this method can be applied.   

Also notice that as a consequence that $u\in\mathcal{C}_1(\Omega\times[0,\infty))$ from the hypothesis, meanwhile the expression remain valid in terms of $u$ does imply that at least $\mathcal{G}\in\mathcal{C}_0(\Omega\times[0,\infty))$ therefore ensure that all the differentials required for $\mathcal{G}$ does exist.

From this is logical integrate the total differential over time as suggest the Picard-Lindel\"of Theorem for Existence and uniqueness referencing the following.

\begin{equation*}
u(t,\vec{x})=u^\circ(\vec{x})+\int_{t_0}^t\mathcal{G}(s,\vec{x},u(s,\vec{x}))ds
\end{equation*}

Now is important to measure how does evolve $u$ inside the Ball centered at $(t_0,\vec{x}_0)$ with radius $\delta$, therefore.

\begin{align*}
|u(t,\vec{x})-u(t_0,\vec{x}_0)|&=\left|u^\circ(\vec{x})-u^\circ(\vec{x}_0)+\int_{t_0}^t\mathcal{G}(s,\vec{x},u(s,\vec{x})) ds\right|\\
&\leq\left|u^\circ(\vec{x})-u^\circ(\vec{x}_0)\right|+\left|\int_{t_0}^t\mathcal{G}(s,\vec{x},u(s,\vec{x})) ds\right|\\
&\leq\left|u^\circ(\vec{x})-u^\circ(\vec{x}_0)\right|+\int_{t_0}^t\left|\mathcal{G}(s,\vec{x},u(s,\vec{x}))\right| ds
\end{align*}

Because $u^\circ\in\mathcal{C}_1(\Omega)$, that way is ensure to be Lipschitz. Also is important to notice that as a consequence of the continuity of $\mathcal{G}$ is constrained to the set $I\times u(I)$, then exist a supreme to call allowing the following inequality.

\begin{align*}
|u(t,\vec{x})-u(t_0,\vec{x}_0)|&\leq\left|u^\circ(\vec{x})-u^\circ(\vec{x}_0)\right|+\int_{t_0}^t\left|\mathcal{G}(s,\vec{x},u(s,\vec{x}))\right| ds\\
&\leq A\left|\vec{x}-\vec{x}_0\right|+M(t-t_0)\\
&\leq A\left(|\vec{x}-\vec{x}_0|+\frac{M}{A}|t-t_0|\right)
\end{align*}

Where the Lipschitz constant is $A$ , do not confuse them with the coefficient of the PDE $\vec{A}$, and the supreme related to $\mathcal{G}$ on $I\times u(I)$ is denoted by $M$, notice that this last expression is a metric inside this domain $I$, therefore is valid to write the following.

\begin{equation*}
|u(t,\vec{x})-u(t_0,\vec{x}_0)|\leq A|(t,\vec{x})-(t_0,\vec{x}_0)|
\end{equation*}

Proving this way that $u(t,\vec{x})$ is Lipschitz over $I$, also notice that by setting $|(t,\vec{x})-(t_0,\vec{x}_0)|<\delta=\epsilon/A$ establish that the function is well define for any $\epsilon>0$, proving this way the existence of the solution.

Now to prove the uniqueness of the function, is necessary to prove that for the following operator.

\begin{align*}
u_{n+1}(t,\vec{x})=T(u_n(t,\vec{x}))=u^\circ(\vec{x})+\int_{t_0}^t\mathcal{G}(s,\vec{x},u_n(s,\vec{x}))ds& &u_0(t,\vec{x})=u^\circ(\vec{x})
\end{align*}

Is a contraction on $I$, building up a local but unique solution for the PDE, that can be expanded to the full domain of the function.

Therefore is necessary to prove that this sequence build upon the operator $T$ converge uniformly toward a fixed function, describe by the Banach Fixed Point Theorem, this is equivalent to prove that this sequence is Cauchy convergent.

Therefore analyzing two first contiguous elements of this sequence is obtained.

\begin{align*}
|u_{1}(t,\vec{x})-u_0(t,\vec{x})|&\leq\left|\int_{t_0}^t\mathcal{G}(s,\vec{x},u^\circ(\vec{x})) ds\right|\\
&\leq\int_{t_0}^t\left|\mathcal{G}(s,\vec{x},u^\circ(\vec{x}))\right| ds\\
&\leq M|t-t_0|\\
&\leq A|(t,\vec{x}_0)-(t_0,\vec{x}_0)|\\
&\leq A|(t,\vec{x})-(t_0,\vec{x}_0)|\\
&< A\delta\\
&<\epsilon
\end{align*}

For $|(t,\vec{x})-(t_0,\vec{x}_0)|<\delta=\epsilon/A$ given $\epsilon>0$. Now by comparing another arbitrary contiguous element of the sequence.

\begin{align*}
|u_{n+1}(t,\vec{x})-u_n(t,\vec{x})|&=\left|\int_{t_0}^t\left[\mathcal{G}(s,\vec{x},u_n(s,\vec{x}))-\mathcal{G}(s,\vec{x},u_{n-1}(s,\vec{x}))\right] ds\right|\\
&\leq\int_{t_0}^t\left|\mathcal{G}(s,\vec{x},u_n(s,\vec{x}))-\mathcal{G}(s,\vec{x},u_{n-1}(s,\vec{x}))\right| ds
\end{align*}

Now as a consequence that $\mathcal{G}\in\mathcal{C}_1(I\times u(I))$, $\mathcal{G}$ is a Lipchitz function on $I\times u(I)$, done this assumption is implied the following.

\begin{align*}
|u_{n+1}(t,\vec{x})-u_n(t,\vec{x})|&\leq L|t-t_0||u_n(t,\vec{x}))-u_{n-1}(t,\vec{x}))|
\end{align*}

Where $L$ is the Lipchitz constant, therefore by iterating this result until reach the first difference of this sequence.

\begin{align*}
|u_{n+1}(t,\vec{x})-u_n(t,\vec{x})|&\leq L^n|t-t_0|^n|u_1(t,\vec{x}))-u_{0}(t,\vec{x}))|\\
&\leq AL^n|t-t_0|^n|(t,\vec{x})-(t_0,\vec{x}_0)|\\
&\leq \frac{A^{n+1}L^n}{M^n}|(t,\vec{x})-(t_0,\vec{x}_0)|^{n+1}\\
&\leq \frac{A^{n+1}L^n}{M^n}\delta^{n+1}\\
&<\epsilon
\end{align*}

For any $\epsilon>0$ with this particular $\delta=\min\{\epsilon,(M/L)\}/A$ ensure that the difference between 2 consecutive elements of the sequence are well define and contained on $u(I)$.

Now taking the difference between any element of the sequence where $m,n\in\mathbb{N}$.

\begin{align*}
|u_m(t,\vec{x})-u_n(t,\vec{x})|&=\left|\sum_{i=n}^{m-1}u_{i+1}(t,\vec{x})-u_{i}(t,\vec{x})\right|\\
&\leq\sum_{i=n}^{m-1}\left|u_{i+1}(t,\vec{x})-u_{i}(t,\vec{x})\right|\\
&\leq|u_1(t,\vec{x}))-u_{0}(t,\vec{x}))|\sum_{i=n}^{m-1} (L|t-t_0|)^i\\
&\leq A\delta\sum_{i=n}^{m-1}\left(\frac{LA}{M}|(t,\vec{x}_0)-(t_0,\vec{x}_0)|\right)^i\\
&\leq A\delta\sum_{i=n}^{m-1}\left(\frac{LA}{M}|(t,\vec{x})-(t_0,\vec{x}_0)|\right)^i\\
&< A\delta\sum_{i=n}^{m-1}\left(\frac{LA\delta}{M}\right)^i
\end{align*}

Now let $\delta=min\{\epsilon,M/L\}/2A$, therefore for all $\epsilon>0$ imply the following.

\begin{align*}
|u_m(t,\vec{x})-u_n(t,\vec{x})|< A\delta\sum_{i=n}^{m-1}\left(\frac{1}{2}\right)^i= A\delta\left(\frac{1-\left(\frac{1}{2}\right)^m}{1-\frac{1}{2}}\right)<2A\delta=\epsilon
\end{align*}

Therefore concluding this proof, we have extended the Picard-Lindel\"of Theorem to PDE's with this specific and particular form describe in equation (\ref{GPDE}), thanks to the Banach fixed point theorem by proving this sequence of functions converge over Cauchy criteria to a specific function that solves the PDE.

As an example let's recall the heat equation, according to the analysis done in that problem from equation (\ref{Heat-Grad}) and (\ref{Heat-temp}) the total derivative becomes.

\begin{align*}
D_Tu=u\left(\frac{\alpha}{4}\sum_{i=1}^n\left(\frac{\partial_tf_i}{f_i}D_T(x_i-x_i^\circ)^2+\left(\frac{\partial_tf_i }{f_i}\right)^2(x_i-x_i^\circ)^2\right)+\frac{1}{2}\sum_{i=1}^n\frac{\partial_tf_i}{f_i}\right)
\end{align*}

That when is recall the condition for the heat equation expressed on equation (\ref{ODE-Heat-System}) and it's solution does imply.

\begin{align}
D_Tu=-u\left(\frac{\alpha}{4}\sum_{i=1}^nD_T\left(\frac{(x_i-x_i^\circ)^2}{(t+1)}\right)+\frac{n}{2(t+1)}\right)\label{Heat-Total}
\end{align}

By taking the partial derivative with respect to the function $u$ does hold.

\begin{align*}
\partial_uD_Tu=-\left(\frac{\alpha}{4}\sum_{i=1}^nD_T\left(\frac{(x_i-x_i^\circ)^2}{(t+1)}\right)+\frac{n}{2(t+1)}\right)
\end{align*}

By evaluating this result over the domain we notice that is well define, therefore under this criteria the solution for the Heat equation under the studied conditions has to be unique.

Notice that this can be done as a consequence of the HRM because give expressions for partial differential along the plausible directions embedded on  $\Omega$ implying the behavior of the gradient and the temporal partial differential from the behavior of the PDE, allowing that way a similar conditions to what happen on a ODE.

\section{From Homogeneous to Non-Homogeneous diffusive PDE.}

\hspace{0.2 in} Now to extend this method for the general case, can be use the following substitution $C(t,\vec{x},u(t,\vec{x}))=C^\prime(t,\vec{x},u(t,\vec{x}))-h(t,\vec{x})/u(t,\vec{x})$ on equation (\ref{GPDE}) implying the following.

\begin{equation}
\Delta u+\vec{A}(t,\vec{x},u)\cdot\nabla u(t,\vec{x})+B(t,\vec{x},u)\frac{\partial u(t,\vec{x})}{\partial t}+\left(C^\prime(t,\vec{x},u(t,\vec{x}))-\frac{h(t,\vec{x})}{u(t,\vec{x})}\right)u=0\label{NHGPE}
\end{equation}

Notice that this last expression is equivalent to.

\begin{align*}
\Delta u+\vec{A}(t,\vec{x},u)\cdot\nabla u+B(t,\vec{x},u)\frac{\partial u}{\partial t}+C^\prime(t,\vec{x},u)u=h(t,\vec{x})
\end{align*}

Therefore compatible with the application of the theorem previously studied on this paper, implying as a consequence of equation (\ref{reduce}).

\begin{equation}
u\left(-\vec{A_\Gamma}\cdot\nabla\Gamma-B_\Gamma\frac{\partial\Gamma}{\partial t}\right)+\Gamma\left((\vec{A}-\vec{A_\Gamma})\cdot\nabla u+(B-B_\Gamma)\frac{\partial u}{\partial t}+(C^\prime-C_\Gamma^\prime)u\right)+(\Gamma-1)h=2\nabla u\cdot\nabla\Gamma\label{Non-Homogeneous-Reduce}
\end{equation}

Notice from this last expression is added a term referencing the $h$ function, allowing this way a direct extension for the HRM for a non linear case, regardless if it is an homogeneous PDE or not being extremly useful for the analysis of this particular type of dynamical system.

Now as a last example to apply all the potential for the method will be studied the three dimensional case of the Navier-Stokes Equations for a incomprehensible fluid with an arbitrary pressure function that depends on time and space and a null force, for a smooth initial data that has a null divergence.

This Equations are characterized to be almost impossible to analyze with normal methods that leads to a many dead ends along their studies becoming a Millennium Prize Problems as a consequence of it difficulty, therefore becomes a perfect target to study the extension of the HRM.

Now to start this last example will be necessary to describe the equations that will be work assuming that the density of the fluid is 1 and the viscosity coefficient $\nu>0$ for the directions $i\in\{1,2,3\}$ orthogonal between each other.

\begin{align}
\frac{\partial u_i}{\partial t}+u\cdot\nabla u_i&=\nu\Delta u_i-\frac{\partial p}{\partial x_i}\label{NSEQ-i}\\
\nabla\cdot u&=0\label{NSEQ-NDIV}
\end{align}

Where the following border condition holds for $i\in\{1,2,3\}$.

\begin{align}
\nabla u_i\cdot\hat{r}=0 \quad \forall \vec{r}\in\partial\Omega\label{NS-Border Condition}
\end{align}

And it's initial data is constrained by the following condition for some $K\in\mathbb{R}^+$.

\begin{align}
\left|\left|u^\circ(\vec{r})\right|\right|\leq\frac{C}{(1+r)^K}\label{NS-InitialData-Constrained}
\end{align}

Notice that as a consequence of the divergence condition from equation (\ref{NSEQ-NDIV}), does hold the following regardless of the reference frame.

\begin{align*}
0=\int_Su\cdot\hat{r}dS
\end{align*}

Therefore $u(t,\vec{r})\cdot\hat{r}=0$ for all element on the border for any $t\geq 0$, implying the following, regardless of the border condition for the gradient of $u_i$. Implying that for this particular case where the gradient of $u_i$ is null for along the radial direction for all $i$ along the border of the domain, Allows that by calculating the volume integral over the domain does hold as a consequence of the divergence theorem and the null divergence of $u$.

\begin{align*}
\int_V\left(\partial_tu_i+\partial_ip\right)dV=\int_V\left(\nu\Delta u_i-u\cdot\nabla u_i\right)dV=\int_S\left(\nu\nabla u_i-u_iu\right)\cdot \hat{n}dS=0
\end{align*} 

This can only holds if and only if.

\begin{align*}
\partial_tu_i+\partial_ip=0
\end{align*} 

This result is independent from the index $i$, therefore holds.

\begin{align}
\partial_tu+\nabla p=0\label{U-temporal-partial}
\end{align}

By calculating the inner product with respect the function $u$.

\begin{align*}
u\cdot\partial_tu+u\cdot\nabla p=0
\end{align*}

That is equivalent to.

\begin{align*}
\partial_t||u||^2=-2\nabla\cdot(pu)
\end{align*}

Notice that by computing the volume integral, as a consequence of the divergence theorem and the border condition establish by the incomprehensibility of the fluid does hold.

\begin{align*}
\int_V\partial_t||u||^2dV=0
\end{align*}

Therefore implying that the kinetic energy of the fluid remain constant over time, therefore the amount of energy of the system over the domain, is constrained by the initial data describe by the inequality (\ref{NS-InitialData-Constrained}) implying the following.

\begin{align*}
\left|\left|u(t,\vec{r})\right|\right|=\left|\left|u^\circ(\vec{r})\right|\right|\leq\frac{C}{(1+r)^K}
\end{align*}

Therefore the total amount of energy is described by the following.

\begin{align*}
\int_V\left|\left|u\right|\right|^2dV\leq&\int_V\left(\frac{C}{(1+r)^K}\right)^2dV=4C^2\pi\int_0^R\left(\frac{r}{(1+r)^K}\right)^2dr\\
=&4C^2\pi\int_0^R\left(\frac{1}{(1+r)^{K-1}}-\frac{1}{(1+r)^{K}}\right)^2dr\\
=&\frac{4C^2\pi}{(1+R)^{2K}}\left(\frac{(1+R)^3}{3-2K}+\frac{(1+R)}{1-2K}-\frac{(1+R)^2}{1-K}\right)-\frac{4C^2\pi}{(3-2K)(1-2K)(1-K)}
\end{align*}

Therefore the energy is finite over the domain of the function for any $K>3/2$, Also notice that as a consequence of such $K$ does hold, that the energy is bounded by.

\begin{align*}
\int_V\left|\left|u\right|\right|^2dV<\frac{4C^2\pi}{(2K-3)(2K-1)(K-1)}
\end{align*}

Also does hold the following from this analysis on the original PDE.

\begin{align}
\nu\Delta u_i-u\cdot\nabla u_i=0\label{NS-Disect}
\end{align}

This allows the existence of a free divergence vector potential $\vec{B}_i$ for each direction, such that holds the following as a consequence that $u$ has null divergence.

\begin{align}
\nu\nabla u_i=u_iu+\vec{B_i}\label{gradient-NS}
\end{align}

Therefore the jacobian matrix of $u$ denoted by $\mathbf{J}(u)$ is the following by concatenating the respective directions building up the next $3\times3$ matrix.

\begin{align*}
\nu\mathbf{J}(u)=u\otimes u+B\label{Jacobian-U}
\end{align*}

Where $B$ is a matrix encapsulating the behavior of the potential vectors, where it components are describe as $B_{ij}=\vec{B}_i\cdot\hat{e_j}$. 

Therefore let's impose that $u$ has the next form using a potential scalar function $\phi$.

\begin{align}
u=-\nu\frac{\nabla\phi}{\phi}
\end{align}

This lead to the following PDE set.

\begin{align*}
\nu\mathbf{J}(u)=-\nu^2\mathbf{J}\left(\frac{\nabla\phi}{\phi}\right)=-\nu^2\left(\frac{\mathbf{H}(\phi)}{\phi}-\frac{\nabla\phi\otimes\nabla\phi}{\phi^2}\right)=-\nu^2\frac{\mathbf{H}(\phi)}{\phi}+u\otimes u
\end{align*}

Where $\mathbf{H}(\phi)$ is the hessian matrix of the potential function $\phi$, Therefore.

\begin{align}
\frac{\mathbf{H}(\phi)}{\phi}=-\frac{B}{\nu^2}
\end{align}

let's notice that this expression can be multiply by the transpose of the gradient of the function $\phi$, leading the following expression.

\begin{align*}
(\nabla\phi)^T\mathbf{H}(\phi)=-\left(\frac{\phi}{\nu^2}\right)(\nabla\phi)^TB
\end{align*}

Notice that this expression are vectors on $\mathbb{R}^3$, but not just that this expression is equivalent to.

\begin{align}
\nabla(||\nabla\phi||^2)^T=-\frac{(\nabla\phi^2)^TB}{\nu^2}\label{Grad-GradPHYNorm}
\end{align} 

Notice that as a consequence that $u$ has to be a free divergence field, does imply from equation (39) the the following statement holds.

\begin{align*}
\frac{||\nabla\phi||^2}{\phi^2}=\frac{\Delta\phi}{\phi}=\frac{tr(\mathbf{H}(\phi))}{\phi}=-\frac{tr(B)}{\nu^2}
\end{align*}

Therefore implying.

\begin{align}
||\nabla\phi||^2=-\frac{\phi^2tr(B)}{\nu^2}\label{Norm-grad-phi}
\end{align}

Using this result on the equation (\ref{Grad-GradPHYNorm}) leads to the following expression.

\begin{align}
\left(\frac{\nabla\phi}{\phi}\right)^T+\frac{1}{2}\left(\frac{\nabla tr(B)}{tr(B)}\right)^T=\left(\frac{\nabla\phi}{\phi}\right)^T\frac{B}{tr(B)}\label{phi-PDE}
\end{align} 

Therefore implying the following equivalence.

\begin{align*}
\left(\left(tr(B)I-B\right)^T\right)\left(\frac{\nabla\phi}{\phi}\right)+\frac{1}{2}\nabla tr(B)=\vec{0}
\end{align*}

If assume $\mathbf{det}(tr(B)I-B)\neq0$, therefore, this partiular matrix has an inverse and also it transpose therefore.

\begin{align*}
\left(\frac{\nabla\phi}{\phi}\right)=-\frac{1}{2}\left(\left(tr(B)I-B\right)^T\right)^{-1}\nabla tr(B)
\end{align*}

Therefore the candidate for $u$ becomes.

\begin{align*}
u=\frac{\nu}{2}\left(\left(tr(B)I-B\right)^T\right)^{-1}\nabla tr(B)
\end{align*}

Equivalent to.

\begin{align*}
u=\frac{\nu}{2}\left(\left(tr(B)I-B\right)^{-1}\right)^T\nabla tr(B)
\end{align*}

Where the associated inverse matrix of our interest correspond to be.

\begin{align*}
(tr(B)I-B)^{-1}=\frac{(-1)^{m+n}\left(M_{mn}(tr(B)I-B)\hat{e}_m\otimes\hat{e}_n\right)^T}{\mathbf{det}(tr(B)I-B)}
\end{align*}

Where $M_{mn}$ refer to the minor matrix determinant associated to the designated coordinate by excluding the $m$-th row and the $n$-th column of the original matrix. Because of the complexity to read the behavior of the $B$ matrix because of the behavior of the potential vectors embedded on it. Therefore would be useful to change the notation system to Einstein notation. Implying that the vector function $u$ does become.

\begin{align*}
u(t,\vec{r})=\frac{\nu(-1)^{m+n}M_{mn}(tr(B)I-B)\partial_ntr(B)}{2\mathbf{det}(tr(B)I-B)}\hat{e}_m
\end{align*}

Notice that $m$ is a dummy index therefore by setting $m=i$ does hold.

\begin{align*}
u(t,\vec{r})=\frac{\nu(-1)^{i+n}M_{in}(tr(B)I-B)\partial_ntr(B)}{2\mathbf{det}(tr(B)I-B)}\hat{e}_i
\end{align*}

Notice that as a consequence of the indexes there exist 9 elements to take on consideration in particular 3 for each component of the vector function $u$, that comes from the $n$ index. Therefore the next step would be computing each element of the function vector $u$ in terms of the elements of the matrix $B$.

Therefore for this function to be well define requires that the matrix $tr(B)I-B$ does have a non zero determinant.

By computing the determinant do we have the following.

\begin{align*}
\mathbf{det}(tr(B)I-B)=\frac{\varepsilon^{ijk}\varepsilon^{i^\prime j^\prime k^\prime}}{6}\left(tr(B)\delta_i^{i^\prime}-B_{ii^\prime}\right)\left(tr(B)\delta_j^{j^\prime}-B_{jj^\prime}\right)\left(tr(B)\delta_k^{k^\prime}-B_{kk^\prime}\right)
\end{align*}

This expression after some algebra is equivalent to.

\begin{align*}
\mathbf{det}(tr(B)I-B)=tr(B)tr(\mathbf{C}(B))-det(B)
\end{align*}

Where $\mathbf{C}(B)$ indicate the co-factor matrix associated to the $B$ matrix. Therefore the only thing left to determine is the minor determinant $M_{in}(tr(B)I-B)$

\begin{align*}
M_{in}(tr(B)I-B)=\frac{\varepsilon^{ijk}\varepsilon^{nj^\prime k^\prime}}{2}\left(tr(B)\delta_j^{j^\prime}-B_{jj^\prime}\right)\left(tr(B)\delta_k^{k^\prime}-B_{kk^\prime}\right)
\end{align*}

Implying the following.

\begin{align*}
M_{in}(tr(B)I-B)=M_{in}(B)+tr(B)\left(B_{ii}\delta_n^i-B_{ji}\delta_n^j-B_{ki}\delta_n^k\right)
\end{align*}

Therefore $u$ becomes.

\begin{align}
u(t,\vec{r})=\frac{(-1)^{i+n}\nu}{2}\left(\frac{M_{in}(B)+tr(B)\left(B_{ii}\delta_n^i-B_{ji}\delta_n^j-B_{ki}\delta_n^k\right)}{tr(B)tr(\mathbf{C}(B))-det(B)}\right)\partial_ntr(B)\hat{e}_i\label{NS-Algebraic Sol}
\end{align}

Remember that this expression is contracted thanks to the Einstein notation. Also notice that by setting the time parameter $t=0$ does imply that we have a specific expression that describe the initial data $u^{\circ}(r)$.

Now the HRM would be use to find an alternative construction for the partial derivatives by reducing this expression from a second order PDE to a first order PDE.

Therefore by applying the HRM, with $\vec{A}=-u/\nu$, $B=-1/\nu$ and $C=0$ for a non homogeneous case does imply the following, after some algebra on the original equation describe in equation (\ref{NSEQ-i}).

\begin{align*}
u_i\left(\left(u+u_i(\Gamma-1)\hat{e_i}\right)\cdot\nabla\Gamma+\frac{\partial\Gamma}{\partial t}\right)-\Gamma(\Gamma-1)u_i\hat{e_i}\cdot\nabla u_i+(\Gamma-1)\frac{\partial p}{\partial x_i}=2\nu\nabla u_i\cdot\nabla\Gamma
\end{align*}

By selecting the following harmonic function valid in spherical coordinates $\Gamma=1+f_i/r$ does hold after some algebra reveals over the radial direction how the PDE behave.

\begin{align*}
u_i\left(\left(u+\frac{u_if_i\hat{r}}{r}\right)\cdot\left(\frac{-\hat{r}}{r}\right)+\frac{1}{f_i}\frac{\partial f_i}{\partial t}\right)+\frac{\partial p}{\partial x_i}=\nabla u_i\cdot\left(\left(1+\frac{f_i}{r}\right)u_i\hat{e_i}-\left(\frac{2\nu\hat{r}}{r}\right)\right)
\end{align*}
 
Because the gradient pressure depend on the temporal partial derivative of $u$ as shown in equation (\ref{U-temporal-partial}), therefore this result becomes.

\begin{align}
\partial_tu_i=u_i\left(\left(u+\frac{u_if_i\hat{r}}{r}\right)\cdot\left(\frac{-\hat{r}}{r}\right)+\frac{\partial_tf_i}{f_i}\right)-\nabla u_i\cdot\left(\left(1+\frac{f_i}{r}\right)u_i\hat{e_i}-\left(\frac{2\nu\hat{r}}{r}\right)\right)\label{reduce-NS}
\end{align}

Therefore this can be resume on the following.

\begin{align*}
\nabla u_i&=\frac{u_iu+\vec{B_i}}{\nu}
\end{align*}
\begin{align*}
\partial_tu_i&=u_i\left(\left(u+\frac{u_if_i\hat{r}}{r}\right)\cdot\left(\frac{-\hat{r}}{r}\right)+\frac{\partial_tf_i}{f_i}\right)-\nabla u_i\cdot\left(\left(1+\frac{f_i}{r}\right)u_i\hat{e_i}-\left(\frac{2\nu\hat{r}}{r}\right)\right)
\end{align*}

Notice that from equation (\ref{NS-Algebraic Sol}) $B_{ij}\in C^1(\mathbb{R}^3\times\left[0,\infty\right))$ for all indexes $(i,j)$, because the function $u\in C^1(\mathbb{R}^3\times\left[0,\infty\right))$ by hypothesis, because each component of this function $u$ belong to $C^1(\mathbb{R}^3\times\left[0,\infty\right))$, therefore by calculating the respective derivatives on the gradient and the temporal partial derivative to ensure they belong to $C^1(\mathbb{R}^3\times\left[0,\infty\right)-\vec{0}\times[0,\infty))$, this would lead to the conclusion that the function $u\in C^2(\mathbb{R}^3\times\left[0,\infty\right)-\vec{0}\times[0,\infty))$, therefore, $B_{ij}\in C^2(\mathbb{R}^3\times\left[0,\infty\right)-\vec{0}\times[0,\infty))$, by repeating this procedure is concluded that the function $u\in C^\infty(\mathbb{R}^3\times\left[0,\infty\right)-\vec{0}\times[0,\infty))$, and because this result does holds, does imply that $p\in C^\infty(\mathbb{R}^3\times\left[0,\infty\right)-\vec{0}\times[0,\infty))$.

Notice that because is just one troublesome element in the domain for a fixed time $t$. For that specific time $t$ at the vicinity of this particular point all elements of the function $u$ are well posed, therefore by calculating the limit over that vicinity of $\vec{0}$ representing the source of the flow leads to a fixed value, by iterating this result over all $t$ the function $u$ collapse to a temporal function for each component, because of that $u$ has to remain smooth over the vicinity of the domain for all $t$, leads to the conclusion that this temporal functions explicit for each coordinate of the function $u$ has to be smooth over time, this last result is true thanks to the energy boundary condition, otherwise if that bound didn't exist we cannot tell if this solution can be corrected or not because the amount of energy over this particular point could be infinite leading to nonsense results.

Therefore is proven that $u,p\in C^\infty(\mathbb{R}^3\times\left[0,\infty\right))$. 

Also the total derivative of the Navier-Stokes equations become.

\begin{align*}
D_tu_i=\frac{\nu^2\vec{B_i}\cdot u-u_itr(B)}{\nu^3}+\frac{(u_iu+2\vec{B_i})\cdot\hat{r}}{r}-\frac{u_i^2f_i}{r^2}+\frac{u_i\partial_tf_i}{f_i}-\frac{u_i^3+u_iB_{ii}}{\nu}\left(1+\frac{f_i}{r}\right)
\end{align*}

Notice that $f_i$ has to be at least a $C^2(\left[0,\infty\right))$ function because by extracting the gradient from equation (\ref{reduce-NS}) thanks to equation (\ref{U-temporal-partial}) the gradient of pressure along that direction and the partial time derivative are interchangeable therefore allowing the construction of the Laplacian in terms of $u$ and $p$ leaving a partial solution of this equation implicating a second derivative for the function $f_i$.

So far we have shown that this function $u$ and $p$ are smooth functions on the domain, the energy of the system is bounded in the domain, therefore finite and we have build an specific algebraic form that depend on the potential divergence free vectors that were added when the gradient of each coordinate where build, therefore the only thing remain to proof is the uniqueness of the solution.

For this will be used the proof by contradiction.

Lets assume that exist $u^1$ and $u^2$ such that are solution to the Navier Stokes equations and has the same initial data and the border data embedded but $u^1\neq u^2$.

Let's define $\omega:=u^2-u^1$ by calculating the energy associated to $\omega$.

\begin{align*}
E_\omega(t)=\int_V\left|\left|\omega(t,\vec{r})\right|\right|^2dV
\end{align*}

By taking the time derivative of this  quantity does imply.

\begin{align*}
\partial_tE_\omega(t)=&\int_V\omega\cdot\partial_t\omega dV=\int_V\omega\cdot\partial_t\left(u^2-u^1\right)dV=\int_V\left(u^2-u^1\right)\cdot\nabla\left(p^1-p^2\right)dV\\
=&\int_V\left(u^2\cdot\nabla p^2-u^1\cdot\nabla p^2\right)dV-\int_V\left(u^2\cdot\nabla p^1-u^1\cdot\nabla p^1\right)dV\\
=&\int_V\nabla\cdot\left(u^2p^2-u^1p^2\right)dV-\int_V\nabla\cdot\left(u^2p^1-u^1p^1\right)dV\\
=&\int_Sp^2\left(u^2-u^1\right)\cdot\hat{n}dS-\int_Sp^1\left(u^2-u^1\right)\cdot\hat{n}dS=0
\end{align*}

Therefore notice that $E_\omega$ is a fixed constant over time, thanks to the border condition that comes from the free divergence field for the function $u$, therefore $u^1$ and $u^2$ should also hold as solution of the Navier-Stokes equations, therefore does hold.

\begin{align*}
E_\omega(t)=E_\omega(0)=\int_V\left|\left|\omega(0,\vec{r})\right|\right|^2dV=\int_V\left|\left|u^2(0,\vec{r})-u^1(0,\vec{r})\right|\right|^2dV=\int_V\left|\left|u^\circ(\vec{r})-u^\circ(\vec{r})\right|\right|^2dV=0
\end{align*}

Therefore $\omega$ has to be a zero vector, implying that $u^1=u^2$, showing this way that under the assumption of a second solution for the Navier-Stokes equation exist a contradiction, thus this result regardless that we don't have an explicit equation to describe the Navier-Stokes Equations, because we need to define the behavior of the elements of the matrix $B$, we have shown that has a smooth and unique solution therefore well posed under the specified conditions.

But this does not provided evidence for the uniqueness of $p$, that is the main reason that on the previous prof did appear $p^1$ and $p^2$.

To prove that the function $p$ that satisfies the Navier-Stokes equations is unique can be prove with ease assuming that $p$ is not unique, therefore exist $p^1$ and $p^2$ that satisfy the Navier-Stokes equation.

\begin{align*}
\nabla(p^2-p^1)=\nabla p^2-\nabla p^1=-\partial_tu+\partial_tu=0
\end{align*}

Notice that this expression establish that 

\begin{align*}
p^2(t,\vec{r})=p^1(t,\vec{r})+f(t)
\end{align*}

As a consequence that they have to coincide for all the elements of the border of the domain it imply that $f$ has to be the zero function for all $t\geq0$, allowing the conclusion that $p^1$ and $p^2$ are the same function, therefore unique.

To express this solution notice that we have equation (\ref{NS-Algebraic Sol}) that is an algebraic expression of the particular solution for the data we have, therefore the only thing left is to determine each component of the matrix $B$ described in the development of this, taking in consideration the amount of data as well the complexity of following each permutation, I do recommend the use of computer assistance to develop more on this matter to adjust the matrix $B$ using Helmholtz vector decomposition to describe the behavior of turbulence as well.

Notice fist that I've been integrating over the domain $V$ and $S$ over the course of this demonstrations to maintain readability, but this result can be easily change by $V=\mathbb{R}^3$ and $S=\partial V$ representing this way the whole space.

Also notice that for a compressible fluid this equation has to be modified in terms of how does behave the density factor over space and time variable.                                                                                                                                                                                                                                                                                                                                                                                                                                                                                            

In this case the HRM give us an alternative view on the PDE to determine it's uniqueness, but because the total derivative is not define at the source of the flux we can not conclude nothing over this point by the extension of the Picard-Lindel\"of, therefore necessary the use of the energy method, right now so many subjects are dependent on the Navier-Stokes equation and by this analysis this topics can be studied with more ease, being the main reason I did some examples referring on nonlinear PDE's instead of focusing on the linear cases.

\section{Some thoughts on the process developing HRM}

\hspace{0.2in} As you might read the HRM has a lot of potential not just to reach new horizons to understand dynamical systems, but also on the development of new and outstanding methods to solve other types of problems concerning to PDE analysis and subtopics that depend on this subject like physical and chemical sciences or ecological systems.

Another thing to remark it might be hard to digest this type of method at first, by the simple reason that it seams that we are deriving a "non-related" function or set of functions to integrate over the spatial coordinates of the domain of the functions, but this results come from a geometrical resources such as the divergence theorem and the Green's identities, being the main reason that the metric and the inner product are necessary that exist with in the domain of the function.

Also is important to acknowledge that most of the methods are related to the Sobolev's spaces and weak derivatives that I try to avoid at all cost, not because they were useless, because they are not, but because, I was thinking on providing a more straight forward approach for diffusive dynamical systems, similar to Laplace's Transform for ODE's because of how important are this type of models on modern sciences.

Next will be presented an extensive bibliography, not because I wanted, but because was necessary, if I'm proposing a new way to develop and analyze PDE's trough the HRM, I should need and extensive bibliography to make sure under human constrain that my works remain original, and also relevant to the mathematical community. concluding this way this 5 year research I've made.

\bibliographystyle{authordate1}
\bibliography{HRM_NSAS}

\nocite{habermanPDE1}
\nocite{asmarPDE2}
\nocite{jostPDE3}
\nocite{moonPDE4}
\nocite{coltonPDE5}
\nocite{elsgoltsPDE6}
\nocite{ortizPDE7}
\nocite{zillPDE8}
\nocite{miklavcicPDE9}
\nocite{boycePDE10}
\nocite{St1PDE11}
\nocite{St2PDE12}
\nocite{St3PDE13}
\nocite{shubinPDE14}
\nocite{shubinPDE15}
\nocite{pratapPDE16}
\nocite{diBenedettoPDE17}
\nocite{cheviakovPDE18}
\nocite{struwePDE19}
\nocite{maleshkoPDE20}
\nocite{ottoPDE21}
\nocite{sauvignyPDE22}
\nocite{bergmanPDE23}
\nocite{powersPDE24}
\nocite{fuscoPDE25}
\nocite{rassiasPDE26}
\nocite{taylorPDE27}
\nocite{egorovPDE28}
\nocite{fedoryukPDE29}
\nocite{agranovichPDE30}
\nocite{shubinPDE31}
\nocite{dunikovPDE32}
\nocite{agranovichPDE33}
\nocite{egorovPDE34}
\nocite{rosingerPDE35}
\nocite{polyaninPDE36}
\nocite{woolseyPDE37}

\nocite{bworldODE1}
\nocite{burtonODE2}
\nocite{perkoODE3}
\nocite{JordanODE4}

\nocite{batemanMPA1}

\nocite{ladyzhenskayaNS1}
\nocite{fernandezNS2}
\nocite{landauNS3}
\nocite{lukaszewiczNS4}

\nocite{garrityManifolds1}
\nocite{mendelsonManifolds2}
\nocite{doCarmoManifolds3}
\nocite{ishamManifods4}
\nocite{ivorraManifolds5}
\nocite{marsdenManifolds6}
\nocite{stoneManifolds7}

\nocite{lovelockTensor1}
\nocite{bishopTensor2}
\nocite{spencerTensor3}
\nocite{borisenkoTensor4}

\nocite{segreDiff1}
\nocite{griffithsDiff2}

\nocite{hwangAnalisis1}
\nocite{taoAnalisis2}
\nocite{fernandezAnalisis3}
\nocite{courantAnalisis4}
\nocite{courantAnalisis5}

\end{document}